\documentclass[11pt]{amsart}
\usepackage{amsfonts,amsmath,amssymb}
\usepackage[all]{xy}
\usepackage{hyperref}
\begin{document}

\newcommand{\ba}{{\bf a}}

\newcommand{\bb}{{\bf b}}
\newcommand{\bpi}{{\boldsymbol{\pi}}}
\newcommand{\oh}{{\mathfrak o}}
\newcommand{\m}{{\mathfrak m}}
\newcommand{\jnf}{{\rm inf}}
\newcommand{\C}{{\mathbb C}}
\newcommand{\F}{{\mathbb F}}
\newcommand{\Gal}{{\rm Gal}}
\newcommand{\Q}{{\mathbb Q}}
\newcommand{\T}{{\mathbb T}}
\newcommand{\Sp}{{\rm Sp}}
\newcommand{\W}{{\sf W}}
\newcommand{\sK}{{\sf K}}
\newcommand{\Z}{{\mathbb Z}} 
\newcommand{\End}{{\rm End}}

\newcommand{\LTp}{{{\rm LT}(\pi)}}
\newcommand{\ltp}{{\sf{lt}(\bpi)}}
\newcommand{\MU}{{{\rm M}{\mathbb U}}}
\newcommand{\con}{{\rm Con}}
\newcommand{\N}{{\mathcal N}}

\parindent=0pt
\parskip=6pt

\newcommand{\ie}{\textit{ie}\,}
\newcommand{\eg}{\textit{eg}\,}
\newcommand{\se}{{\sf e}}
\newcommand{\cf}{{\textit{cf}\,}}

\title{Notes on $\delta$-algebras and prisms in homotopy theory}

\author[J Morava]{J Morava}

\address{Department of Mathematics, The Johns Hopkins University,
Baltimore, Maryland}

\begin{abstract}{J McClure's Dyer-Lashof operation \cite{16} in $p$-adic $K$-theory defines, in particular, a prismatic structure on the complex representation ring of the circle group. Work of  Ando, Rezk, Stapleton, and others generalizes this to define a canonical lift of Frobenius for the Lubin-Tate spectra of \cite{13}. We suggest that recent work of K Ito and S Marks on $\oh_L$ -  typical (or $L$ -adic) prisms may extend this to the $\oh_L$-adic spectra $\sK(\oh_L)$ of \cite{19}.} \end{abstract}

\maketitle 

{\bf \S \;  Local conventions} following \cite{15}(Ch 6) roughly: $L$ is a field of $p$-adic numbers, with 
\[
\infty > n = [L:\Q_p], L \supset \oh_L \supset \m_L = (\pi), \oh_L/\m_L = k_L, |k_L| = q = p^f, f | n .
\]

If $\pi_L = \pi \in \oh_L$ thus parameterizes $L$, then there is an associated Lubin-Tate formal group law $\LTp/\oh_L$ with
\[
X +_\pi Y \in \oh_L[[X,T]], \oh_L \ni a \to [a]_\pi(T) \in \End_{\oh_L}(\LTp)
\]
having mod $\pi$ reduction $\ltp/k_L$ with
\[
X +_\bpi Y \in k_L[[X,T]], \oh_L \ni a \to [a]_\bpi (T) \in  \End_{k_L}(\ltp) ,
\]
the formal $\oh_L$-module (\ie with complex multiplication) structure map in the first case being an isomorphism and injective in the second. Under reduction mod $\pi$,  $[\pi]_\pi(T) \mapsto \bpi$, embedding $L$ in a canonical division $\Q_p$-algebra of isogenies. \bigskip
 
 \; $\bullet$ {\bf Spectral background} Quillen, together with Lazard, associates to a (one-dimensional) formal group law over a commutative ring $A$, a Hirzebruch multiplicative genus of complex oriented manifolds, \ie a homomorphism from the complex cobordism ring $\MU^*$ to $A[v,v^{-1}]$. In particular, the middle arrow in the factorization of the homomorphism 
 \[
\chi_L(\pi) :  \MU^* \cong {\rm Lazard}^* \to E^*_{\ltp,k_L} \otimes_{W(k_L)} \oh_L \to  \oh_L^*[v^{\pm 1}]
\]
(\ie $\C P(q^k - 1) \mapsto \pi^{-k}q^k$, otherwise $0$, classifying $\LTp/\oh_L$) collapses the fundamental constructions of Lubin and Tate in algebra with the subsequent work of mathematical generations \cite{20} in homotopy theory, culminating in the existence of complex-oriented $E_\infty$ ring spectra such as 
\[
\pi_0 E^*_{\ltp,k_L} \cong W(k_L)[[\dots]][v^{\pm 1}] 
\]
with a formal group law canonically identified with the universal deformation ring of, for example, $\ltp/k_L$. 

Charles Rezk \cite{22}(\S 11.13) considers power operations in tensor (two-)categories of quasicoherent sheaves over a category of deformations of a one-dimensional formal group law (eg $\ltp/k_L$ together with its Frobenius) regarded as a category of comodules over a bialgebra $\Gamma$ (of semi-linear operations with respect to a monad {\bf T} on the category of $E_*$-modules), \cf \cite{6}. His generalization (Theorem B, \S 1.4) of Clarence Wilkerson's criterion \cite{24} identifies \cite{21}(\S 1.12, 14.3) the Frobenius endomorphism of a large class of $E_n$-algebras as a Hecke operator $\sigma \in \Gamma \otimes \F_p$ (\cf Mathew Ando \cite{1}, Nora Ganter \cite{8}(\S 3.8)), generalizing McClure's Dyer-Lashof work in Atiyah's \cite{2} $p$-adic context. Note that these are {\bf un}stable Hecke operations on cohomology algebras of {\bf spaces}, rather than spectra. 

In \cite{21}(Prop 3.6) Nathaniel Stapleton interprets this as a canonical lift of Frobenius for Goerss-Hopkins-Miller-Rezk Lubin-Tate theories associated to one-dimensional formal groups over finite fields, as described above: in particular, providing it with a $\delta_{\Q_p}$ structure as described in \S \2 below.\bigskip

\; $\bullet$ {\bf Genuine Morava $\sK(\oh_L)$-theories} In \cite{19}(\S 2.3.2 Prop), Baas-Sullivan Koszul resolutions \cite{7} kill the kernel of $\chi_L(\pi)$ to define small-batch artisanal $E_2$ ringspectra $\sK(\oh_L)$. The succeeding section \S 3.1 identifies 
\[
\sK(\oh_L)_* \sK(\oh_L)  \cong  \oh_L \otimes_{E^*_{\ltp,k_L}} {\rm HAut}_{\oh_L}(\LTp) \otimes_{\oh_L} \Lambda^*(\m_L \otimes \Z/pZ)^\vee
\]
as a $\mu_L$ - graded exterior $\oh_L$-algebra extension of Andrew Baker's \cite{3}(Prop 2.5) $\oh_L$-adic Banach Hopf algebra
\[
{\rm HAut}_{\oh_L}(\LTp) \cong \con(\oh_L^\times,\oh_L)   \cong    H{\mu_L} \otimes_{\oh_L} \con(\oh_L,\oh_L),
\]
where
\[
\con(\oh_L,\oh_L) \supset \N(\oh_L) := \oh_L[\theta_i \:|\: i \geq 0] 
\]
is generated by the continuous numerical polynomial functions
\[
\theta_k(T) =  - \pi^{-k}[\sum_{0 \leq i \leq k-1} \pi^i \theta_i(T)^{q^{k-i}} - T]
\]
(\eg $\theta_1(T) = \pi^{-1}(T - T^q)$) from $\oh_L$ to itself. Evidently $\N(\oh_L)$ admits the identity map
\[
\phi_q(T) = T^q + \pi \cdot \theta_1(T) = T
\]
as a lift of the $q$-Frobenius on the quotient ring $\con(\oh_L,k_L)$, making it into a $\delta_L$-algebra. \bigskip

\newpage

{\bf \S \;  an $L$ - adic prism}, following A Baker, K Ito, S Marks and N Stapleton\bigskip

\; $\bullet$ {\bf Proposition} {\it If $A \in$ local $\oh_L$ - algebras, then its (commutative) ring $\W(A)$ of length two ramified Witt vectors is a topologically nilpotent thickening
\[
\xymatrix{
0 \ar[r] & A_\pi \ar[r]^-{a \to (a,0)} & \W(A) = A \times A \ar[r]^-{\ba \to a_1} & A \ar[r] & 0}
\]
of $A$}, in coordinates $\ba = (a_0,a_1)$, where $A_\pi$ is the nonunital ring with $A$ as underlying abelian group, but with product 
\[
A_\pi \ni a,b \mapsto \pi \dot ab .
\] 

Following Kazuhiro Ito \cite{14}(lemma 2.2.6) and Samuel Marks \cite{17}(\S 2.3.1), Spec $\W$ is the ($\oh_L$-algebra) - scheme defined by composition laws
\[
\ba +_\W \bb = (a_0 + b_0, a_1 + b_1  + \pi^{-1} \langle a_0,b_0 \rangle_q ) 
\]
\[
\ba \times_\W \bb = (a_0 b_0, \pi a_1 b_1 + (a_1,b_1)^{\rm T} \cdot (a_0^q ,b_0^q))
\]
(where T denotes transpose and $\langle x,y \rangle_q = x^q + y^q - (x + y)^q$).\bigskip

Similarly following James Borger \cite{4}(Theorem A, Fig 1), Sp $\W(A)$ is the pushout of two copies of $\Sp A$ along the Frobenius endomorphism of $A/\pi A$, presented as
\[
\Sp A/\pi A \Rightarrow \Sp_{\oh_L} A \times_{\Sp_{\oh_L}} \Sp_{\oh_L} A \to \W(A) .
\]
An $\oh_L$-module homomorphism 
\[
A \ni a \mapsto (a,\delta_L(a)) \in \W(A) : \Sp_{\oh_L} \W(A)  \to \Sp_{\oh_L} A
\]
splitting the exact sequence defining the ring $\W$ is a section of a principal $A_\pi^\times$ `bundle' over $\Sp A$, defining something analogous to a normal family of deformations for $A$, cf \cite{9}(\S 1) re B Dwork's lemma. Borger notes [\S 1.19] that $\delta_L$ -algebra structures avoid the choice implicit in a `lift' of Frobenius. \bigskip

\; $\bullet$ {\bf Example} The ($\pi$-typical) Lubin-Tate construction \cite{17}(\S 4.1) defines a ring homomorphism  
\[
T \mapsto \phi_q(T) = [\pi]_L(T) = T^q + \pi  \delta_L(T) : \oh_L[[T]] \to \oh_L[[T]]] 
\]
or, equivalently, a $\delta_L$-structure
\[
\delta_L(T) = \pi^{-1}[[\pi]_L(T) - T^q]
\] 
such that $\phi_q(T) \equiv T^q$ mod $\pi$. The image $\bpi$ of  $[\pi]_L (T) \in k_L[[T]]$ is not just a ring endomorphism of $k_L[[T]]$, but an endomorphism of a formal group of finite height. As an element of the $\Q_p$-division algebra of isogenies of $\ltp$, its minimal polynomial $E(\bpi)$ identifies $\oh_L  \cong W(k_L)[\pi]/E(\pi)$ as a $W(k_L)$-algebra.\bigskip

\; $\bullet$ {\bf Definition} as in Marks, \cf Ito \cite{14}(lemma 2.3.8): If
\[
q_n(T) = [\pi^{n-1}](T)^{-1} \cdot [\pi^n](T) \equiv T^{(q-1)q^{n-1}} + \dots \in \oh_L[[T]] 
\]
then $\oh_L[[T]] \ni q_1(T) = T^{-1} [\pi](T) \equiv \pi$ mod $T$. The $(\pi,T)$-adic completeness of $\oh_L[[T]]$ then implies the $(\pi,q_n(T))$-adic completeness of $A$ with respect to the ideal $(q_n(T))$ \dots

Indeed
\[
\pi = q_1(T) + T \cdot f(T)
\]
with $f(T) \in \oh_L[[T]]$, so after applying $\phi^n$ we have 
\[
\pi  = q_{n+1}(T) + [\pi^n](T) \cdot f([\pi^n](T)) =   
\]
\[
=  q_1(T) + q_n(T) [\pi^{n-1}](T) \cdot  f([\pi^n](T)) \in (q_n(T),q_{n+1}(T)) \subset \oh_L[[T]] .
\]

That is, if $I = (q_n(T))$, then $\pi \in I + \phi(I) \cdot A$ as in \cite{17}(Def 2.3.1, lemma 4.1).\bigskip

\; $\bullet$ {\bf Proposition} {\it With a Lubin-Tate group law as above, $\forall (n \geq 1)$ the pair $(A = \oh_L[[T]],(q_n(T)))$ is an ($\oh_L$-typical) prism}. \bigskip

The following Prop 4.10 defines an interesting perfection 
\[
(\oh_L[[T]],(q_n(T)) \to (A_\jnf(\oh_{L^\infty}),\ker \theta) \dots .
\]
\bigskip

\; $\bullet$ The original version of this note attempted to raise the question of existence of lift of Frobenius (\ie a $\delta_L$-structure) for  $\sK(\oh_L)$ as a consequence of Rezk's work, extending Stapleton in the unramifield case: a geometrically natural lift of Frobenius on height one $E_*$-theory suggests interpreting $\sK(\oh_L)$ as taking values\begin{footnote}{including a Dennis \cite{18} trace $\sK(\oh_L) \to {\rm THH}(\oh_{L^\infty}:\Z_p)$, $L^\infty$ being a maximal totally ramified abelian extension of $L$ or maybe its completion}\end{footnote} in Adams/Artin/Atiyah/Quillen/Iwasawa 
\[
\oh_L[[\Gal(L^\infty/L)]] \cong \oh_L[[\oh_L^\times]]
\]
algebras, presumably related to Fontaine's period ring $A_\jnf(\oh_{L^\infty})$ - with Frobenius lifts (\ie a $\delta_L$-algebra structure).

{\bf Conjecture} $\oh_L \mapsto \sK(\oh_L)$ is the result of some more natural construction, perhaps exhilarating an interesting Gal$(L/\Q_p)$ - action, when applicable.\bigskip

\; $\bullet$ deepest {\bf Thanks} to Andy Baker, Denem Manam, and Nat Stapleton for helpful correspondence and conversation. Hallucinatory howlers are the author's entire responsibility.\bigskip

\; $\bullet$ Executive Summary : Monkey find prism.

\bibliographystyle{amsplain}

 \end{document}